\documentclass{gtart_h}

\def\ifplaintex{\expandafter\ifx\csname documentclass\endcsname\relax}

\def\gtp{{\mathsurround=0pt\it $\cal G\mskip-2mu$eometry \&\ 
$\cal T\!\!$opology $\cal P\!$ublications}}  

\def\recd{{\small Received:\qua\receiveddate\ifx\reviseddate\relax
\else\qquad Revised:\qua\reviseddate\fi\par}} 

\def\lognumber#1{\def\thelognumber{#1}}
\def\volumenumber#1{\def\thevolumenumber{#1}}
\def\volumeyear#1{\def\thevolumeyear{#1}}
\def\papernumber#1{\def\thepapernumber{#1}}
\def\pagenumbers#1#2{\def\startpage{#1}\def\finishpage{#2}}
\def\published#1{\def\publishdate{#1}}

\def\received#1{\def\receiveddate{#1}}

\def\accepted#1{\def\accepteddate{#1}}
\def\asciititle#1{\def\theasciititle{#1}}

\def\asciiauthors#1{\def\theasciiauthors{#1}}
\def\asciiaddress#1{\def\theasciiaddress{#1}}
\def\asciiemail#1{\def\theasciiemail{#1}}

\def\coverauthors#1{\def\thecoverauthors{#1}}
\long\def\asciiabstract#1{\long\def\theasciiabstract{#1}}
\def\asciikeywords#1{\def\theasciikeywords{#1}}

\let\\\par\let\thelognumber\relax\let\thevolumenumber\relax
\let\thepapernumber\relax\let\thevolumeyear\relax\let\startpage\relax
\let\finishpage\relax\let\publishdate\relax\let\receiveddate\relax
\let\reviseddate\relax\let\accepteddate\relax\let\theasciititle\relax
\let\theasciiauthors\relax\let\theasciiaddress\relax
\let\theasciiabstract\relax\let\theasciikeywords\relax

\let\thecoverauthors\relax\let\theasciiemail\relax

\ifplaintex
\font\logobig=cmssbx10 scaled 3836
\font\logomed=cmssbx10 scaled 2557
\else
\font\logobig=cmssbx10 scaled 4200
\font\logomed=cmssbx10 scaled 2800
\fi

\long\def\makeagttitle{   
\count0=\startpage
\agt\hfill      
\hbox to 45truept{\vbox to 0pt{\vglue -13truept{\logomed A\kern -.37em{\logobig 
T}\kern -.38em G}\vss}\hss}
\break
{\small Volume \thevolumenumber\ (\thevolumeyear)
\startpage--\finishpage\nl
Published: \publishdate}

\vglue .25truein

{\parskip=0pt\leftskip 0pt plus
1fil\def\\{\par\smallskip}{\Large\bf\thetitle}\par\medskip} \vglue
0.05truein

{\parskip=0pt\leftskip 0pt plus 1fil\def\\{\par}{\sc\theauthors}
\par\medskip}%
 
\vglue 0.03truein 

{\small\leftskip 25truept\rightskip 25truept{\bf Abstract}\stdspace\theabstract

{\bf AMS Classification}\stdspace\theprimaryclass
\ifx\thesecondaryclass\relax\else; \thesecondaryclass\fi\par
{\bf Keywords}\stdspace \thekeywords\par}\vglue 7truept

}   

\ifplaintex
\font\phead=cmsl9 scaled 950
\font\pnum=cmbx10 scaled 913
\font\pfoot=cmsl9 scaled 950
\headline{\vbox to 0pt{\vskip -4.5mm\line{\small\phead\ifnum
\count0=\startpage ISSN 1472-2739 (on-line) 1472-2747 (printed)
\hfill {\pnum\folio}\else\ifodd\count0\def\\{ }%
\ifx\theshorttitle\relax\thetitle\else\theshorttitle\fi\hfill{\pnum\folio}
\else\def\\{ and }{\pnum\folio}\hfill\ifx\theshortauthors\relax\theauthors
\else\theshortauthors\fi\fi\fi}\vss}}
\footline{\vbox to 0pt{\vglue 0mm\line{\small\pfoot\ifnum\count0=\startpage
\copyright\ \gtp\hfill\else
\agt, Volume \thevolumenumber\ (\thevolumeyear)\hfill\fi}\vss}}
\else
\font=cmsl9 scaled 1050
\font\lnum=cmbx10 
\font=cmsl9 scaled 1050
\makeatletter
\def\@oddhead{{\small\lhead\ifnum\count0=\startpage ISSN 1472-2739 
(on-line) 1472-2747 (printed)\hfill {\lnum\number\count0}\else\ifodd\count0
\def\\{ }\ifx\theshorttitle\relax \thetitle \else\theshorttitle\fi\hfill
{\lnum\number\count0}\else\def\\{ and }{\lnum\number\count0}
\hfill\ifx\theshortauthors\relax 
\theauthors\else\theshortauthors\fi\fi\fi}}\def\@evenhead{\@oddhead}
\def\@oddfoot{\small\lfoot\ifnum\count0=\startpage\copyright\ \gtp\hfill\else
\agt, Volume \thevolumenumber\ (\thevolumeyear)\hfill\fi}
\def\@evenfoot{\@oddfoot}
\makeatother
\fi
\let\maketitlepage\makeagttitle
\let\makeshorttitle\maketitlepage
\let\maketitle\maketitlepage

\newwrite\gtoutfile
\long\gdef\makeheadfile{  
{\def\\{, }\def\s{ }
\immediate\openout\gtoutfile head.xxx
\immediate\write\gtoutfile{Proxy-for: \ifx\theasciiauthors\relax
\theauthors\else\theasciiauthors\fi\s<\ifx\theasciiemail\relax\theemail\else\theasciiemail\fi>}
\immediate\write\gtoutfile{\noexpand\\}
\immediate\write\gtoutfile{Authors: \ifx\theasciiauthors\relax
\theauthors\else\theasciiauthors\fi}
{\def\\{ }\immediate\write\gtoutfile{Title: \ifx\theasciititle\relax
\thetitle\else\theasciititle\fi}}
\immediate\write\gtoutfile{Subj-class: GT or SG, GR etc}
\immediate\write\gtoutfile{MSC-class: \theprimaryclass\ifx\thesecondaryclass\relax\else, \thesecondaryclass\fi}
\immediate\write\gtoutfile{Journal-ref: Algebraic and Geometric Topology \thevolumenumber\s
(\thevolumeyear) \startpage-\finishpage}
\immediate\write\gtoutfile{Comments: Published by Algebraic and
Geometric Topology at}
\immediate\write\gtoutfile{\s\s\s  http://www.maths.warwick.ac.uk/agt/AGTVol\thevolumenumber/agt-\thevolumenumber-\thepapernumber.abs.html}
\immediate\write\gtoutfile{\noexpand\\}
\immediate\write\gtoutfile{}
\ifx\theasciiabstract\relax
\immediate\write\gtoutfile{\theabstract}\else
\immediate\write\gtoutfile{\theasciiabstract}\fi
\immediate\write\gtoutfile{}
\immediate\write\gtoutfile{\noexpand\\}
\immediate\write\gtoutfile{}
\immediate\closeout\gtoutfile}}  

\def\maketitlepage{\makeagttitle\makeheadfile}
\let\makeshorttitle\maketitlepage
\let\maketitle\maketitlepage

\lognumber{12}
\volumenumber{4}
\volumeyear{2004}
\papernumber{12}
\published{10 April 2004}
\pagenumbers{199}{217}
\received{6 October 2003}
\accepted{31 March 2004}

\usepackage{amsmath, amssymb, graphicx, graphs, verbatim}

\newtheorem{thm}{Theorem}[section]
\newtheorem{cor}[thm]{Corollary}
\newtheorem{lem}[thm]{Lemma}
\newtheorem{prop}[thm]{Proposition}

\theoremstyle{definition}

\newtheorem{conj}{Conjecture}

\newtheorem{rem}[thm]{Remark}

\numberwithin{equation}{section}

\newcommand{\al}{\alpha}

\newcommand{\Ga}{\Gamma}

\newcommand{\La}{\Lambda}

\newcommand{\om}{\omega}

\newcommand{\Si}{\Sigma}

\newcommand{\s}{\mathbf s}

\renewcommand{\t}{\mathbf t}

\newcommand{\C}{\mathbb C}
\newcommand{\Z}{\mathbb Z}
\newcommand{\N}{\mathbb N}
\newcommand{\Q}{\mathbb Q}
\newcommand{\R}{\mathbb R}

\newcommand{\CP}{{\mathbb C}{\mathbb P}}

\newcommand{\del}{\partial}

\newcommand{\lra}{\longrightarrow}
\newcommand{\hf}{{{\widehat {HF}}}} 
 
\newcommand{\Li}{\mathbb {L}}
\newcommand{\Spin}{{\rm {Spin}}}

\DeclareMathOperator{\PD}{PD}

\begin{document}

\title{Seifert fibered contact three--manifolds via surgery}
\asciititle{Seifert fibered contact three-manifolds via surgery}

\author{Paolo Lisca\\Andr\'{a}s I. Stipsicz}
\coverauthors{Paolo Lisca\\Andr\noexpand\'{a}s I. Stipsicz}
\asciiauthors{Paolo Lisca\\Andras I. Stipsicz}

\address{Dipartimento di Matematica\\
Universit\`a di Pisa \\I-56127 Pisa, Italy } 

\secondaddress{R\'enyi Institute of Mathematics\\
Hungarian Academy of Sciences\\
H-1053 Budapest\\ 
Re\'altanoda utca 13--15, Hungary}
\asciiaddress{Dipartimento di Matematica,
Universita di Pisa\\I-56127 Pisa, Italy\\and\\Renyi 
Institute of Mathematics,
Hungarian Academy of Sciences\\
H-1053 Budapest,
Realtanoda utca 13--15, Hungary}

\gtemail{\mailto{lisca@dm.unipi.it}{\qua\rm 
and\qua}\mailto{stipsicz@math-inst.hu}}
\asciiemail{lisca@dm.unipi.it, stipsicz@math-inst.hu}

\begin{abstract}
  Using contact surgery we define families of contact structures on
  certain Seifert fibered three--manifolds. We prove that all these
  contact structures are tight using contact Ozsv\'ath--Szab\'o
  invariants. We use these examples to show that, given a natural
  number $n$, there exists a Seifert fibered three--manifold carrying
  at least $n$ pairwise non--isomorphic tight, not fillable
  contact structures.  
\end{abstract}
\asciiabstract{%
  Using contact surgery we define families of contact structures on
  certain Seifert fibered three-manifolds. We prove that all these
  contact structures are tight using contact Ozsath-Szabo
  invariants. We use these examples to show that, given a natural
  number n, there exists a Seifert fibered three-manifold carrying
  at least n pairwise non-isomorphic tight, not fillable
  contact structures.}

\asciikeywords{Seifert fibered 3-manifolds, tight, fillable
contact structures, Ozsvath-Szabo invariants}
\maketitle

\section{Introduction and statement of results}\label{s:intro}
The classification problem for tight contact structures on closed oriented
three--manifolds is one of the driving forces in present day contact topology.
Contact surgery along Legendrian links provides a powerful tool for
constructing contact three--manifolds. Tightness of these structures is,
however, hard to prove, unless the structures can be shown to
be~\emph{fillable}, i.e., can be viewed as living on the boundary of a
symplectic four--manifold satisfying appropriate compatibility conditions. The
question whether any tight contact structure is fillable was open for some
time, until the first tight, nonfillable contact three--manifolds were found
by Etnyre and Honda~\cite{EH}, followed by infinitely many such
examples~\cite{LS1, LS2}. The tightness of those examples was proved using a
delicate topological method called~\emph{state traversal} (see~\cite{H2}). In
this paper we prove tightness by applying the Heegaard Floer theory recently
developed by Ozsv\'ath and Szab\'o \cite{OSzF1, OSzF2, OSz6}.  According to
our main result, tight, not fillable contact structures are more common than
one would expect:

\begin{thm}\label{t:main}
  For any $n\in \N$ there is a Seifert fibered $3$--manifold $M_n$
  carrying at least $n$ pairwise non--isomorphic tight, not
  fillable contact structures.
\end{thm}

The construction of the contact structures in Theorem~\ref{t:main}
relies on contact surgery. We verify nonfillability via the
Seiberg--Witten equations, following the approach of~\cite{LS1, LS2}.
In order to state precisely our results we need a little
preparation.

\subsection*{Contact surgery}
In a given contact three--manifold $(Y, \xi)$ a knot $K\subset
(Y,\xi)$ is~\emph{Legendrian} if $K$ is everywhere tangent to
$\xi$. The framing of $K$ naturally induced by $\xi$ is called
the~\emph{contact framing}. Given a Legendrian knot $K$ in a contact
three--manifold $(Y,\xi)$ and a rational number $r\in\Q$ ($r\neq 0$),
one can perform \emph{contact $r$--surgery} along $K$ to obtain a new
contact three--manifold $(Y',\xi')$~\cite{DG1, DG2}. Here $Y'$ is the
three--manifold obtained by smooth $r$--surgery along $K$, where the
surgery coefficient is measured~{\sl with respect to the contact
framing} defined above, not with respect to the framing induced by a
Seifert surface (which, in general, does not exist). The contact
structure $\xi'$ is constructed by extending $\xi$ from the complement
of a standard neighborhood of $K$ to a tight contact structure on the
glued--up solid torus.  If $r\neq 0$ such an extension always exists,
and for $r=\frac 1k$ $(k\in \Z)$ it is unique \cite{H2}.  When $r=-1$
the corresponding contact surgery coincides with Legendrian surgery
along $K$~\cite{el,G,W}.

Below we outline an algorithm for replacing a contact $r$--surgery on a
Legendrian knot $K$ with a sequence of contact $(\pm 1)$--surgeries on a
suitable Legendrian link.  By~\cite[Proposition~3]{DG2}, contact $r$--surgery
along $K\subset (Y, \xi )$ with $r<0$ is equivalent to Legendrian surgery
along a Legendrian link $\Li =\cup_{i=0}^m L_i$ which is determined via the
following simple algorithm by the Legendrian knot $K$ and the contact surgery
coefficient $r$.  The algorithm to obtain $\Li $ is the following. Let
\[
[a_0+1, \ldots ,a_m],\quad a_0,\ldots a_m\leq -2 
\]
be the continued fraction expansion of $r$. To obtain the first component
$L_0$, push off $K$ using the contact framing and stabilize it $-a_0-2$ times.
Then, push off $L_0$ and stabilize it $-a_1-2$ times.  Repeat the above scheme
for each of the remaining pivots of the continued fraction expansion.  Since
there are $-a_i-1$ inequivalent ways to stabilize a Legendrian knot $-a_i-2$
times, this construction yields $\Pi_{i=0}^m (-a_i-1)$ potentially different
contact structures. According to~\cite[Proposition~7]{DG2}, a contact
$r=\frac{p}{q}$--surgery ($p,q\in \N$) on a Legendrian knot $K$ is equivalent
to a contact $\frac{1}{k}$--surgery on $K$ followed by a contact
$\frac{p}{q-kp}$--surgery on a Legendrian pushoff of $K$ for any integer $k\in
\N$ such that $q-kp<0$.  Therefore, the latter surgery can be turned into a
sequence of Legendrian surgeries, as described above.
By~\cite[Proposition~9]{DG1}, a contact $\frac{1}{k}$--surgery ($k\in \N$) on
a Legendrian knot $K$ can be replaced by $k$ contact $(+1)$--surgeries on $k$
Legendrian pushoffs of $K$. 

In conclusion, any contact rational $r$--surgery ($r\neq 0$) can be replaced
by contact $(\pm 1)$--surgery along a Legendrian link (which is not
necessarily uniquely specified); for a related
discussion see also~\cite{DGS}.

\subsection*{Statement of results}
In the following, we shall denote by
\[
M(g,n;(\alpha _1 ,\beta _1), \ldots , (\alpha _k ,\beta _k))
\]
the Seifert fibered $3$--manifold obtained by performing
$(-\frac{\alpha _1 }{\beta _1})$--, $\ldots, (-\frac{\alpha _k }{\beta
  _k})$--surgeries along $k$ fibers of the circle bundle $Y_{g,n}\to
\Sigma _g$ over the genus--$g$ surface $\Sigma _g$ with Euler number
$e(Y_{g,n})=n$.  The Seifert invariants
\[
(g,n; (\alpha _1, \beta _1 ), \ldots , (\alpha _k , \beta _k))
\]
are said to be in~\emph{normal form} if 
\[
\alpha_i > \beta _i \geq 1,\quad i=1,\ldots , k. 
\]
Using Rolfsen twists (hence changing $n$ if necessary), any tuple
\[
(g,n; (\alpha_1, \beta_1 ), \ldots , (\alpha_k , \beta_k))
\]
can be transformed into normal form.

\begin{figure}[ht]
\begin{center}
\setlength{\unitlength}{1mm}
\begin{picture}(90, 120)
\put(0,0){\includegraphics[height=12cm]{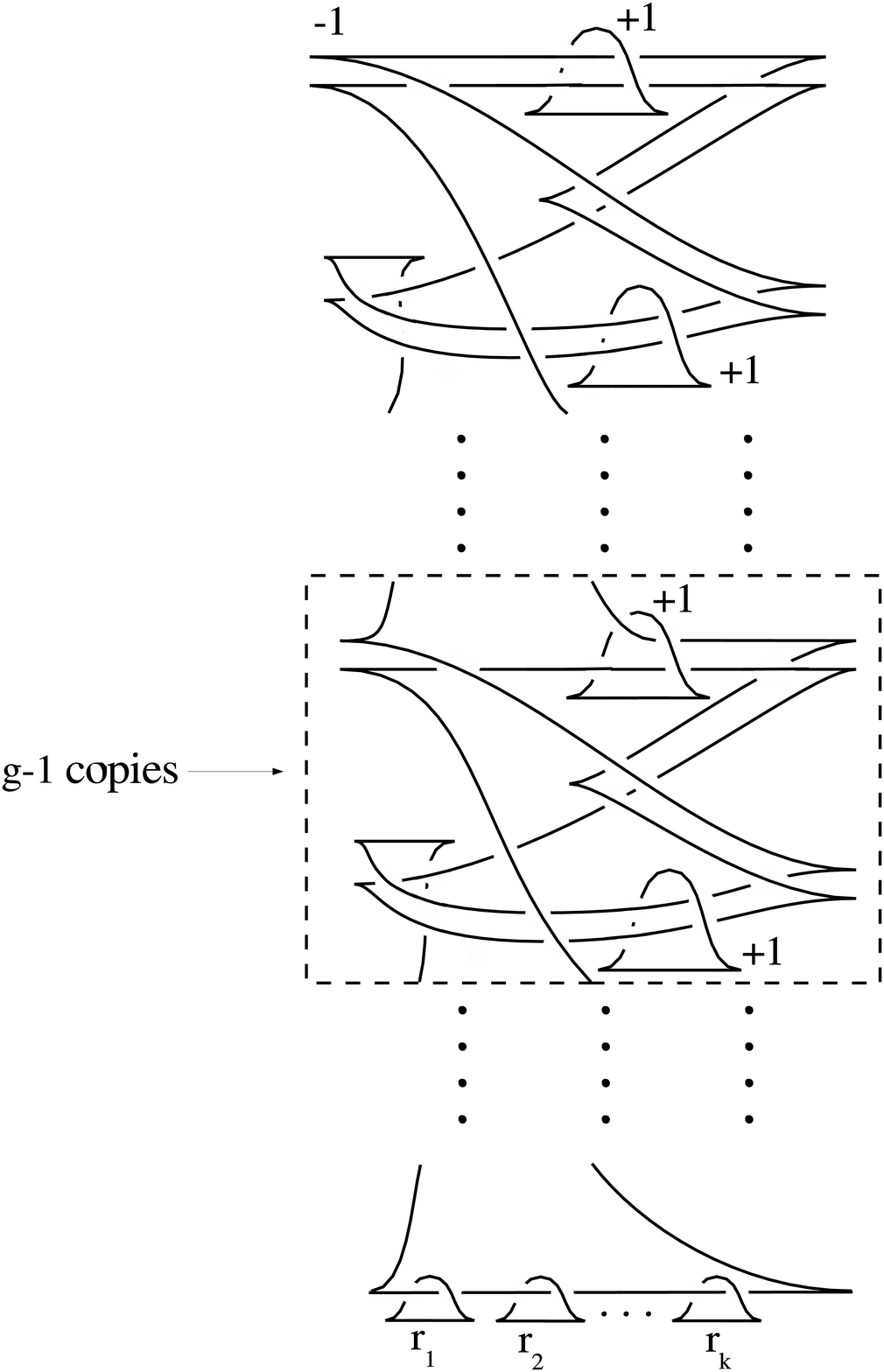}}
\end{picture}    
\end{center}
\caption{Contact structures on Seifert fibered $3$--manifolds}
\label{f:general}
\end{figure}

Consider the family of contact $3$--manifolds defined by the contact
surgery diagrams of Figure~\ref{f:general} (the box is repeated 
$(g-1)$--times, $g\geq 1$).

Throughout the paper we shall assume: 
\begin{equation}\label{e:conditions}
g\geq 1, \quad \frac 12 \leq r_1<1,\quad 
r_i<0, \quad i=2,\ldots ,k \quad (r_i\in \Q). 
\end{equation}
Under the assumptions~\eqref{e:conditions} one can write the
coefficients as:
\begin{equation}\label{e:coefficients}
r_1=\frac{(n-2g+1)\alpha _1+\beta _1}{(n-2g+2)\alpha _1+\beta _1}, 
\quad r_i=\frac{\beta _i -\alpha _i}{\beta _i},
\end{equation}
where 
\[
n\geq 2g,\quad 
\alpha _1> \beta _1\geq 0, \quad 
\alpha _i > \beta _i \geq 1,\ i=2,\ldots, k.
\]
Converting the contact surgery coefficients into smooth coefficients, after
$(n-2g+1)$ Rolfsen twists on the $r_1$--framed unknot we conclude that the
$3$--manifolds underlying the contact structures given by
Figure~\ref{f:general} are of the form:
\begin{equation}\label{e:seifert}
M(g,n;(\alpha _1 ,\beta _1), \ldots , (\alpha _k ,\beta _k)),
\quad n\geq 2g. 
\end{equation}
Moreover, if $\beta _1>0$ the Seifert invariants are in normal form.  Observe
that for $\beta _1=0$ the $(-\frac{\alpha _1}{\beta _1})$--surgery is trivial.

Conversely, given a Seifert fibered $3$--manifold $M$ as
in~\eqref{e:seifert}, Figure~\ref{f:general} provides a contact
structure on $M$ as long as the coefficients $r_i$ defined
by~\eqref{e:coefficients} satisfy the conditions~\eqref{e:conditions}.

Let $\xi_1,\ldots,\xi_t$ denote the contact structures obtained by turning the
diagrams of Figure~\ref{f:general} into contact $(\pm 1)$--surgeries in all
possible ways according to the algorithm described in the previous subsection.
This paper is devoted to the study of $\xi_1,\ldots,\xi_t$. Using the contact
Ozsv\'ath--Szab\'o invariants~\cite{OSz6} we prove:

\begin{thm}\label{t:tight}
  Fix $k\geq 1$, $g\geq 1$, $\frac 12 \leq r_1<1$ and $r_i<0$ for $i=2, \ldots
  , k$. Then, all the contact structures defined by Figure~\ref{f:general} are
  tight.
\end{thm}

It is unclear from the construction whether the contact structures
$\xi_1,\ldots,\xi_t$ are all distinct up to isotopy. Observe that for
$k=1$ and $r_1=\frac{\alpha +1} {2\alpha +1}$ the $3$--manifold underlying 
Figure~\ref{f:general} is $M(g,2g; (\alpha, 1))$.

\begin{thm}\label{t:distinct}
  Given $g\geq 1$ and $n \in \N$, there is an $\alpha\in\N$ such that
  at least $n$ of the contact structures defined by
  Figure~\ref{f:general} for $k=1$ and $r_1=\frac{\alpha +1} {2\alpha
    +1}$ are pairwise non--isomorphic.
\end{thm}

In fact, a more detailed analysis shows that the contact structures
defined by Figure~\ref{f:general} on $M(g,2g; (\alpha, 1))$ are all
distinct up to isotopy (see Section~\ref{s:nonfill}). This leads us to:

\begin{conj}\label{c:differ}
All the tight contact structures defined by Figure~\ref{f:general} and 
satisfying the assumptions~\eqref{e:conditions} are distinct up to isotopy.
\end{conj}

Recall that a contact $3$--manifold $(Y,\xi)$ is~\emph{symplectically
  fillable}, or simply~\emph{fillable}, if there exists a compact
  symplectic four--manifold $(W, \omega)$ such that (i) $\partial W=Y$
  as oriented manifolds (here $W$ is oriented by $\om\wedge\om$) and
  (ii) $\omega\vert_{\xi}\not=0$ at every point of $Y$. Our next
  result concerns fillability properties of some of the contact
  structures under examination.

\begin{thm}\label{t:nonfill1}
  Fix $\al\in\N$ and $g\geq 1$ such that $d(d+1)\leq 2g\leq
  d(d+2)-1$ for some positive integer $d$.  Then, the tight contact
  structures defined by Figure~\ref{f:general} for $k=1$ and
  $r_1=\frac{\al+1}{2\al+1}$ are not symplectically fillable.
\end{thm}

As we show in Section~\ref{s:nonfill}, there is some evidence
supporting the following:
\begin{conj}\label{c:filling}
  No contact structure defined by Figure~\ref{f:general} and
  satisfying conditions~\eqref{e:conditions} is fillable.
\end{conj}

The above results immediately imply Theorem~\ref{t:main}:

\begin{proof}[Proof of Theorem~\ref{t:main}]
Fix $n\in \N$ and $g=1$. Choose $\alpha\in\N$ such that the statement
of Theorem~\ref{t:distinct} holds. The contact structures $\xi _1,
\ldots , \xi _t$ defined by Figure~\ref{f:general} on $M(1,2; (\alpha
, 1))$ are tight by Theorem~\ref{t:tight} and there are at least $n$
pairwise non--isomorphic among them by Theorem~\ref{t:distinct}. By
Theorem~\ref{t:nonfill1} applied with $d=1$ they are also
not fillable. This concludes the proof.
\end{proof}

Our results seem to suggest (see~Section~\ref{s:nonfill}) that a
Seifert fibered $3$--manifold 
\[
M(g, n ; (\alpha _1, \beta _1), \ldots, (\alpha _k , \beta _k ))
\]
with Seifert invariants in normal form should support a tight, not
fillable contact structure if $n\geq 2g>0$. This should be
contrasted with the result of Gompf~~\cite{G}, who showed that a
Seifert fibered $3$--manifold with base genus $g\geq 1$ always carries
a Stein fillable contact structure.

Section~\ref{s:tight} is devoted to the proof of
Theorem~\ref{t:tight}, while Theorems~\ref{t:distinct} and
\ref{t:nonfill1} will be proved in Section~\ref{s:homotopy}.  In
Section~\ref{s:nonfill} we give further evidence supporting
Conjectures~\ref{c:differ} and \ref{c:filling}.

\rk{Acknowledgements} The first author was partially supported by MURST, and
he is a member of EDGE, Research Training Network HPRN-CT-2000-00101,
supported by The European Human Potential Programme. The authors would like to
thank Peter Ozsv\'ath and Zolt\'an Szab\'o for many useful discussions
regarding their joint work. The second author was partially supported OTKA
T034885 and T037735.

\section{Proof of Theorem~\ref{t:tight}}\label{s:tight}

In a remarkable series of papers \cite{OSzF1, OSzF2, OSzF4, OSz6}
Ozsv\'ath and Szab\'o defined new invariants of many low--dimensional
objects --- including contact structures on closed $3$--manifolds. In
this section we apply these invariants to prove Theorem~\ref{t:tight}.

Heegaard Floer theory associates abelian groups $HF^+(Y,\t)$ and $\hf
(Y, \t)$ to a closed, oriented $\Spin^c$ $3$--manifold $(Y, \t)$, and
homomorphisms 
\[
F^+_{W, \s}\co HF^+(Y_1, \t _1)\to HF^+(Y_2, \t _2),\quad
{\widehat F}_{W, \s}\co\hf (Y_1, \t _1)\to \hf (Y_2, \t _2)
\]
to a $\Spin^c$ cobordism $(W, \s)$ between two $\Spin^c$
$3$--manifolds $(Y_1, \t _1)$ and\break $(Y_2, \t _2)$. 

Throughout this paper we shall assume that $\Z/2\Z$ coefficients are
being used in the complexes defining the $HF^+$-- and
$\widehat{HF}$--groups.

Let $Y_{g, -2g}$ be a circle bundle over the genus--$g$ surface
$\Sigma_g$ with Euler number $-2g$ $(g\geq 1)$, and let $D_{g, -2g}$
denote the corresponding disk bundle. Since $H^2 (D_{g, -2g}; \Z )$
has no $2$--torsion, each $\Spin^c$ structure on $D_{g, -2g}$ is
uniquely determined by its first Chern class. Let $\s$ be the unique
$\Spin ^c$ structure on $D_{g, -2g}$ with $c_1(\s)=0$, and denote by
$\t$ the restriction of $\s$ to $Y_{g, -2g}$.

\begin{figure}[ht]
\begin{center}
\setlength{\unitlength}{1mm}
\begin{picture}(95,40)
\put(0,0){\includegraphics[height=4cm]{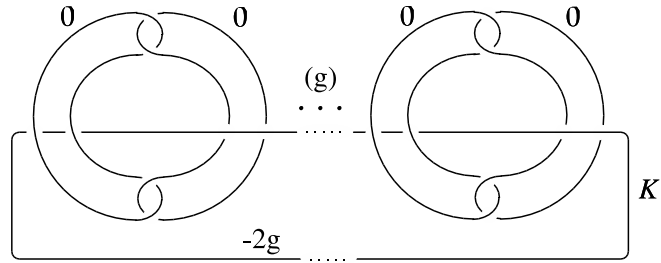}}
\end{picture}    
\end{center}
\caption{The framed knot $K$}
\label{f:cob}
\end{figure}

Let $W$ denote the cobordism from $\#_{2g}(S^1\times S^2)$ to
$Y_{g,-2g}$ given by the attachment of a $4$--dimensional $2$--handle
along the $(-2g)$--framed knot $K\subset \#_{2g}(S^1\times S^2)$ of
Figure~\ref{f:cob}.
Let $\t_0\in\Spin^c (\# _{2g}(S^1\times S^2))$ be the unique $\Spin^c$
structure on $\# _{2g} (S^1\times S^2)$ with vanishing first Chern
class. In~\cite[Lemma~9.17]{abs} it is proved that there is an
isomorphism
\[
HF^+(\#_{2g} (S^1\times S^2),\t_0)\lra HF^+(Y_{g, -2g},\t)
\]
which can be written as a sum of maps $\sum_{\s} F^+_{W,\s}$ over
the set of $\Spin^c$ structures on $W$ which restrict to $\t_0$ and
$\t$. Application of the $5$--lemma to the long exact sequence
connecting $HF^+(Y_{g, -2g}, \t )$ and $\hf (Y_{g, -2g}, \t)$
immediately yields the following:
\begin{lem}\label{l:abs}
The homomorphism 
\[
\sum_{ \{ \s\in \Spin^c(W)\ |\ \s|_{\del W}=(\t_0,\t) \} }{\widehat
  F}_{W,\s}\co \hf (\# _{2g} (S^1\times S^2), \t_0)\to \hf (Y_{g, -2g}, \t ),
\]
is an isomorphism. \qed
\end{lem}

\subsection*{Contact Ozsv\'ath--Szab\'o invariants}
Let $(Y,\xi)$ be a closed contact $3$--manifold oriented by $\xi$, and
let $\t_\xi\in\Spin^c(Y)$ be the $\Spin^c$ structure induced by $\xi$.
In~\cite{OSz6}, Ozsv\'ath and Szab\'o define an invariant
\[
c(Y, \xi)\in \hf (-Y,\t_{\xi })
\]
whose main properties are summarized in the following two theorems.

\begin{thm}\label{t:item1}{\rm \cite{OSz6}}\qua 
If $(Y, \xi )$ is overtwisted, then $c(Y, \xi )=0$. If $(Y, \xi )$ is
Stein fillable then $c(Y, \xi )\neq 0$.  In particular, for the
standard contact structure $(S^3, \xi _{st})$ we have $c(S^3, \xi
_{st})\neq 0$.\qed
\end{thm}

\begin{thm}\label{t:item2}
Suppose that $(Y_2,\xi_2)$ is obtained from $(Y_1,\xi _1)$ by a
contact $(+1)$--surgery. Then we have 
\[
F_{-W} (c(Y_1, \xi_1))= c(Y_2,\xi_2), 
\]
where $-W$ is the cobordism induced by the surgery with reversed
orientation and $F_{-W}$ is the sum of $\sum_\s {\widehat F}_{-W,\s}$ over
all Spin$^c$ structures $\s$ extending the Spin$^c$ structures induced
on $-Y_i$ by $\xi_i$, $i=1,2$.  In particular, if $c(Y_2, \xi_2)\neq
0$ then $(Y_1, \xi_1)$ is tight.
\end{thm}

\begin{proof} Let us assume that we are performing contact
  $(+1)$--surgery along the Legendrian knot $K \subset (Y_1,\xi_1)$.  Then,
  there is an open book decomposition $(F,\phi)$ on $Y_1$ compatible with
  $\xi_1$ in the sense of Giroux and such that $K$ lies on a page. In fact,
  the proof of \cite[Theorem~3]{Gi} shows that the $1$--skeleton of any
  contact cellular decomposition of $(Y_1,\xi_1)$ is contained in a page of a
  compatible open book.  Since $K$ can be assumed to lie in the $1$--skeleton
  of a contact cellular decomposition of $(Y_1,\xi_1)$, the conclusion
  follows.  Moreover, up to refining the decomposition, we may assume that $K$
  is not homotopic to the boundary of the page. Then, an open book for
  $(Y_2,\xi_2)$ is given by $(F,\phi')$, where $\phi'=\phi\circ R^{-1}_K$ and
  $R_K$ is the right--handed Dehn twist along $K$.  The first part of the
  statement now follows applying~\cite[Theorem~4.2]{OSz6}. The second part of
  the statement follows immediately from the fact that the invariant of an
  overtwisted contact structure vanishes.
\end{proof}
Theorem~\ref{t:item2} immediately yields:
\begin{cor}\label{c:surg}
  If $c(Y_2,\xi _2)\neq 0$ and $(Y_1,\xi _1)$ is obtained from
  $(Y_2,\xi _2)$ by Legendrian surgery along a Legendrian knot, then
  $c(Y_1, \xi _1)\neq 0$. In particular, $(Y_1 ,\xi _1)$ is tight.
\end{cor}
\begin{proof} Let $K\subset (Y_2,\xi _2)$ 
  be the Legendrian knot along which the Legendrian surgery is performed. A
  Legendrian pushoff of $K$ gives rise to a Legendrian knot $\widetilde K$ in
  $(Y_1, \xi _1)$.  By~\cite[Proposition~8]{DG1}, contact $(+1)$--surgery on
  $(Y_1,\xi _1)$ along $\widetilde K$ gives back $(Y_2,\xi _2)$. Therefore, by
  Theorem~\ref{t:item2} $c(Y_2,\xi _2)\neq 0$ implies $c(Y_1, \xi _1)\neq 0$.
\end{proof}

\begin{figure}[ht]
\begin{center}
\setlength{\unitlength}{1mm}
\begin{picture}(113,30)
\put(0,0){\includegraphics[height=3cm]{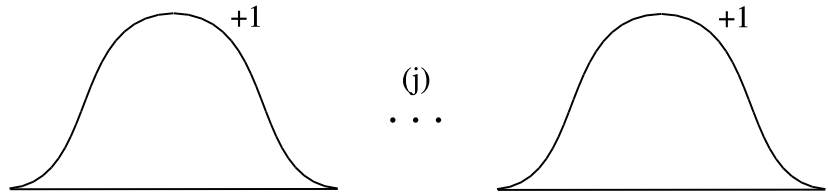}}
\end{picture}    
\end{center}
\caption{The contact $3$--manifold $(Z_j, \eta_j)$}
\label{f:unlink}
\end{figure}

Let $(Z_j, \eta_j)$ be the contact $3$--manifold obtained by
performing contact $(+1)$--surgery on the standard contact
three--sphere along the $j$--component Legendrian unlink depicted in
Figure~\ref{f:unlink}.

\begin{lem}\label{l:unknots}
The contact $3$--manifold $(Z_j, \eta _j)$ given by
Figure~\ref{f:unlink} has non--vanishing contact Ozsv\'ath--Szab\'o
invariant for every $j\geq 0$.
\end{lem}

\begin{proof}
Notice first that $Z_j$ is diffeomorphic to $\# _j (S^1\times S^2)$.  We
will argue by induction on $j$. For $j=0$ we have the standard contact
$3$--sphere, which has non--vanishing contact Ozsv\'ath--Szab\'o invariant by
Theorem~\ref{t:item1}.  Now consider $\eta _{j-1}$ and add the $j$--th
component of the Legendrian unlink to it with contact framing
$(+1)$. Let $-W$ be the corresponding cobordism with reversed
orientation. By~\cite[Theorem~9.16]{OSzF2} the homomorphism
$F_{-W}$ fits into an exact triangle:
\[
\begin{graph}(8,2)
\graphlinecolour{1}\grapharrowtype{2}
\textnode {A}(1,1.5){$\hf (\#_{j-1} (S^1\times S^2))$}
\textnode {B}(7, 1.5){$\hf (\# _j (S^1\times S^2))$}
\textnode {C}(4, 0){$\hf (\#_{j-1} (S^1\times S^2))$}
\diredge {A}{B}[\graphlinecolour{0}]
\diredge {B}{C}[\graphlinecolour{0}]
\diredge {C}{A}[\graphlinecolour{0}]
\freetext (4,1.8){$F_{-W}$}
\end{graph}
\]

In~\cite[Subsection~3.1 and Proposition~6.1]{OSzF2} it is proved
that 
\[
\dim _{\Z/2\Z} \hf (\# _j (S^1\times S^2))=2^j.
\]
Therefore, the exactness of the triangle implies that the map $F_{-W}$
is injective. Since by Theorem~\ref{t:item2} we have 
\[
F_{-W}(c(Z_{j-1}, \eta _{j-1}))=c(Z_j,\eta_j)
\]
and by the inductive assumption $c(Z_{j-1}, \eta _{j-1}))\neq 0$, this
concludes the proof.
\end{proof}

Note that when $k=1$ and $r_1=\frac{1}{2}$, Figure~\ref{f:general}
specifies a unique contact structure $\xi _g$ for every $g$ because
the contact surgery coefficients are of the form $\frac 1k$, $k\in\Z$.
Denote the resulting contact $3$--manifold by $(Y_g, \xi _g)$. It is a
simple exercise to verify that $Y_g$ is an $S^1$--bundle over a
genus--$g$ surface with Euler number $e(Y_g)=2g$.

\begin{prop}\label{p:2g}
  The contact Ozsv\'ath--Szab\'o invariant of $(Y_{g}, \xi _g)$ is nonzero.
\end{prop}

\begin{proof}
  Let $(Y'_g, \xi_g')$ be the contact $3$--manifold given by
  Figure~\ref{f:general} with $k=1$ and $r_1=1$, and perform contact
  $(+1)$--surgery on a pushoff of the $r_1$--framed Legendrian knot $K$.
  According to the algorithm described in Section~\ref{s:intro}, the resulting
  contact structure is $(Y_g, \xi _g)$. Note that $Y_g'$ is diffeomorphic to
  $\#_{2g}(S^1\times S^2)$.  Combining Lemma~\ref{l:unknots} and
  Corollary~\ref{c:surg} we conclude $c(Y_g', \xi_g')\neq 0$. In fact,
  $(Y_g',\xi_g')$ must be the only tight, hence Stein fillable contact
  structure on $\# _{2g}(S^1\times S^2)$. The cobordism given by the handle
  attachment induced by the surgery along $K$ can be easily identified (after
  reversing orientation) with the cobordism appearing in Lemma~\ref{l:abs},
  therefore the non--vanishing of $c(Y_g ' , \xi _g ')$ implies, by
  Theorem~\ref{t:item2}, that $c(Y_{g}, \xi _g)\neq 0$.
\end{proof}

\begin{rem}
  The tightness of the contact structures $\xi_g$ was first proved by
  Honda~\cite{H2} (see also~\cite{LS2}).
\end{rem}

\begin{proof}[Proof of Theorem~\ref{t:tight}]
  Let $K_1, K_2$ denote two Legendrian pushoffs of the $r_1$--framed
  Legendrian unknot $K$ of Figure~\ref{f:general}.  According to the algorithm
  of Section~\ref{s:intro} all contact structures of Figure~\ref{f:general}
  can be given as negative contact surgery on the diagram obtained erasing the
  $r_i$--framed circles $(i=2, \ldots , k)$ from Figure~\ref{f:general} and
  performing contact $(+1)$--surgeries on $K, K_1$ and contact
  $\frac{r_1}{1-2r_1}$--surgery on $K_2$.  (Here we use the assumption $r_i<0$
  for $i=2, \ldots , k$.)  Since $r_1\geq \frac 12$, the surgery coefficient
  of $K_2$ is also negative (or infinity), therefore all the contact
  structures defined by Figure~\ref{f:general} (obeying the restrictions on
  the $r_i$) can be given as Legendrian surgery on $(Y_g,\xi_g)$ for an
  appropriate $g\geq 1$.  Since negative contact surgery can be replaced by a
  sequence of Legendrian surgeries, Corollary~\ref{c:surg} and
  Proposition~\ref{p:2g} imply that these contact structures have
  non--vanishing contact Ozsv\'ath--Szab\'o invariants, hence by
  Theorem~\ref{t:item1} they are tight. This concludes the proof of the
  theorem.
\end{proof}

\section{The proof of non--fillability}
\label{s:homotopy}

Suppose that $(Y, \xi )$ is given by a contact $(\pm 1)$--surgery
diagram and denote the corresponding $4$--manifold by $X$.  Then, the
$\Spin ^c$ structure of the 0--handle of $X$ extends to a $\Spin^c$
structure $\s \in \Spin ^c(X)$ with the property that $\s \vert
_{\partial X}=\t _{\xi}$ and $c_1(\s )$ evaluates on a homology class
$[\Sigma _K] $ given by an oriented surgery curve $K$ as rot$(K)$.
This statement was proved for $(-1)$--surgeries by Gompf~\cite{G} ---
in this case the complex structure of $D^4$ also extends over
the $2$--handles --- and in~\cite{LS2} for the case of
$(+1)$--surgeries; see also \cite{DGS}.

Consider the diagram obtained from Figure~\ref{f:general} for $k=1$
and $r_1= \frac{\alpha +1}{2\alpha +1}$; this diagram represents
contact structures on $M(g, 2g; (\alpha, 1))$. According to the
algorithm outlined in Section~\ref{s:intro}, these contact structures
are also representable by replacing the Legendrian knot $K$ with three
Legendrian pushoffs $K_1, K_2, K_3$ having contact surgery
coefficients $(+1)$, $(+1)$ and $-(\alpha+1)$, respectively.  This
last diagram can be turned into a contact $(\pm 1)$--surgery diagram
by stabilizing the Legendrian curve $K_3$ $\alpha$ times. There are
$(\alpha +1)$ different ways to do this. Choose an orientation for
$K_3$ and define $\xi _r$ as the result of the surgery along the
diagram with rot$(K_3)=r$. (Notice that $r\equiv \alpha$ (mod 2) and
$-\alpha \leq r \leq \alpha$.) The above observation regarding
$\Spin^c$ structures yields:

\begin{lem}\label{l:spinc}
  Let $\s\in \Spin ^c (X)$ be the unique $\Spin^c$ structure such that
  $\langle c_1 (\s), [\Sigma _{K_3}]\rangle =r$ and $\langle c_1 (\s),
  [\Sigma_{j}]\rangle =0$ on the $2$--homology classes defined by the
  remaining surgery circles. Then, the restriction of $\s$ to $\del
  X$ is the $\Spin^c$ structure $\t_{\xi _r}\in \Spin ^c (M(g,2g;
(\alpha , 1)))$ induced by the contact structure $\xi_r$.\qed
\end{lem}

Recall that, since $X$ is simply connected, the Chern class $c_1(\s)$
uniquely specifies the $\Spin ^c $ structure $\s\in \Spin ^c (X)$. For
$M=M(g, 2g; (\alpha , 1))$ let $\mu \in H_1(M; \Z )$ denote the
homology class of the normal circle to the knot $K_3$ --- or,
equivalently, the homology class represented by the singular fiber of
the Seifert fibration. Then, Lemma~\ref{l:spinc} implies that
\[
c_1(\xi_r)=c_1(\t_{\xi_r})=r \PD(\mu). 
\]
In particular, since the order of $\mu $ in $H_1(M; \Z )$ is equal to
$2g\alpha +1$, $\t_{\xi_r}$ is a torsion $\Spin ^c$ structure for all
$r$.

\begin{proof}[Proof of Theorem~\ref{t:distinct}]
  By the classical Dirichlet's theorem on primes in arithmetic
  progressions, there are infinitely many primes of the form $2gm+1$
  as $m$ varies among the natural numbers.  Therefore, we can choose
  natural numbers $a_1,\ldots, a_n$ so that
\[
p_1=2ga_1+1,\ldots,p_n=2ga_n+1
\]
are distinct odd primes. Define $a$ so that 
\[
2ga+1=p_1\cdots p_n. 
\]
If $a$ is odd, let $\al=a$, otherwise let $\al=a(2g+1)+1$. With this
choice $2g\al+1$ is divisible by $p_1\cdots p_n$ and $\al$ is odd.
Therefore, 
\[
\al\equiv p_i\bmod 2,\quad i=1,\ldots,n, 
\]
and we can choose the stabilizations of $K_3$ so that $c_1(\xi _i)=
p_i \mu$.  This implies that the order of $c_1(\xi_i)$ is
$\frac{2g\al+1}{p_i}$, and since the $p_i$'s are all distinct, the
orders of the $c_1(\xi_i)$'s are all different for $i=1, \ldots,
n$. This shows that the contact structures $\xi_i$, $i=1,\ldots, n$,
are pairwise non--isomorphic, concluding the proof.
\end{proof}

The proof of Theorem~\ref{t:nonfill1} will follow the approach used
in~\cite{PLpos} and further exploited in \cite{LS1}. Fix a Seifert
fibration 
\[
M=M(g,n; (\alpha_1, \beta_1 ), \ldots , (\alpha_k , \beta_k))
\to \Sigma _g
\]
over the orbifold $\Sigma _g$. The surface $\Si_g$ can be thought of
as the underlying space of an orbifold with $k$ marked points of
multiplicities $\al_1,\ldots,\al_k$. An orbifold line bundle $L\to
\Sigma _g$ can be pulled back to an honest line bundle ${\overline
  L}\to M$ with torsion first Chern class, and if the invariants
$\alpha _i$ are mutually coprime, all line bundles on $M$ with torsion
first Chern class arise in this way.  An orbifold line bundle $L\to
\Sigma _g$ can be described by its \emph{Seifert data} $(c;\gamma _1,
\ldots , \gamma _k)$, where $c$ is the background degree of $L$ and
the numbers $\gamma _i$ determine the orbifold bundle around the
orbifold points of $\Sigma _g$ (see~\cite[\S 2]{MOY} for further
details).  For example, the orbifold canonical bundle $K_{\Sigma }$
has Seifert data $(2g-2; \alpha _1-1, \ldots , \alpha _k-1)$.  The
\emph{degree} of the orbifold line bundle $L$ is equal by definition
to the rational number
\[
\deg (L)= b+\sum _{i=1}^k \frac{\gamma _i}{\alpha _i}.
\]
For more about Seifert fibered three--manifolds and line bundles on
them see~\cite{MOY, orlik}.

\begin{thm}\label{t:moy}{\rm\cite{MOY}}\qua
The moduli space of Seiberg--Witten solutions for the Seifert
fibered $3$--manifold $M=M(g,2g; (\alpha , 1 ))$ and 
$\Spin ^c$ structure $\t _{\xi _r}
\in \Spin ^c(M)$ 
contains only reducible solutions, for all of which the associated
Dirac operator has trivial kernel. 
\end{thm}
\begin{proof}
  We need to express the $\Spin ^c$ structure $\t_{\xi_r}$ in the
  coordinates used in \cite{MOY} and then appeal to the description of
  the Seiberg--Witten moduli spaces on Seifert fibered $3$--manifolds
  as given in~\cite[Theorem~5.19]{MOY}. In that paper the $\Spin ^c$
  structures are parametrized by their twisting relative to the
  canonical $\Spin ^c$ structure $\t _{\rm can}$ induced by any
  tangent $2$--plane field transverse to the $S^1$--fibration.  As
  explained in~\cite[\S 3]{MOY}, the orbifold disk bundle associated
  to $M$ can be desingularized to a smooth complex surface $X$ with
  $\del X=M$. The group $H_2(X;\Z)$ is generated by the classes of a
  genus--$g$ smooth complex curve $C$ and a smooth rational curve $R$,
  satisfying:
\[
C\cdot C=2g,\quad C\cdot R=1,\quad R\cdot R=-\al.
\]
The restriction to $\del X$ of the complex bundle $TX$ is isomorphic
to the pull--back of
\[
\underline\C\oplus K^{-1}_\Si\to \Si_g,
\]
where $\underline\C$ is the trivial complex line bundle and $K_\Si$ 
is the orbifold canonical bundle of $\Si_g$.

Therefore, denoting by $\s^\C$ the $\Spin ^c$ structure on $X$ induced by
the complex structure, we have $\s^\C\vert_{\partial X} = \t _{\rm can}$ (cf.
text following \cite[Lemma~5.10]{MOY}). The adjunction formula gives:
\[
\langle c_1(X), C\rangle = 2,\quad \langle c_1(X), R\rangle = 2-\al.
\]
Thus, if $\Ga_r\in H^2(X;\Z)$ is a cohomology class satisfying
\[
\langle \Ga_r, C\rangle = -1\quad 
\langle \Ga_r, R\rangle = \frac 12 (r+\alpha -2),
\]
setting $\s_r=\s^\C + \Ga_r$, we have $\s_r|_{\del X}=\t_{\xi_r}$.  
This implies: 
\begin{equation}\label{e:spin-c}
\t_{\xi_r}=\t_{\rm can} + \Ga_r|_{\del X}= 
\t_{\rm can} + \frac 12 (r -\alpha -2)\PD(\mu).
\end{equation}
Now \cite[Theorem~5.19]{MOY} can be restated in the following form, more
convenient for our present purposes. Fix a torsion $\Spin ^c$
structure
\[
\t_k=\t_{\rm can} + k\PD(\mu)\in\Spin ^c (M).
\]
Let $L_k\to \Sigma _g$ be an orbifold line bundle which pulls back to
a line bundle ${\overline L}_k\to M$ with $c_1({\overline
  L_k})=k\PD(\mu)$. Then, the moduli space ${\mathfrak M}_k$ of
Seiberg--Witten solutions on $M$ in the $\Spin ^c$ structure $\t_k$ has a
component of reducible solutions (homeomorphic to the Jacobian torus
of $\Sigma _g$), and by~\cite[Corollary~5.17]{MOY} the associated
Dirac operators have trivial kernels if and only if either $\alpha $
is even or
\begin{equation}\label{e:kernel}
\deg L_k \not\in \frac 12 \deg K_{\Sigma } + 
(2g+\frac{1}{\alpha})\cdot\Z\subset\Q.
\end{equation}\label{e:reduc}
\noindent In addition, ${\mathfrak M}_k$ contains irreducible solutions if and
only if there exists some orbifold line bundle $L\to \Sigma _g$
satisfying:
\begin{equation}\label{e:irred}
\deg L\in [0, \deg K_\Sigma]\setminus\{\frac 12\deg K_\Sigma\},
\quad\deg L\in\deg L_k + (2g+\frac{1}{\alpha })\cdot\Z.
\end{equation}
In view of~\eqref{e:spin-c}, in our case we have:
\[
k=\frac 12 (r -\alpha -2)\equiv 2g\al+ \frac12 (r-\al)\ \ \bmod (2g\al+1).
\]
Therefore, since $r\in [-\al,\al]$, 
\[
\deg K_{\Sigma } = 2g-1-\frac{1}{\alpha}
<\deg L_k=2g+\frac 1{2\al} (r-\al) <
2g+\frac{1}{\alpha }.
\]
It follows that $L_k$ satisfies~\eqref{e:kernel} and there is no
orbifold line bundle $L\to\Sigma _g$ satisfying~\eqref{e:irred}.
Hence, ${\mathfrak {M}}_k$ consists entirely of reducible solutions
with associated Dirac operators having trivial kernels.
\end{proof}

\begin{cor}\label{c:weak}
  Let $(W, \omega )$ be a weak filling of the contact $3$--manifold $(M, \xi
  _r)$. Then, $b_2^+(W)=0$ and the homomorphism $H^2 (W;\R )\to H^2(\partial
  W;\R)$ induced by the inclusion $\del W\subset W$ is the zero map.
\end{cor}
\begin{proof}
  The statement follows from Theorem~\ref{t:moy} in exactly the same
  way as~\cite[Proposition~4.2]{LS1} follows
  from~\cite[Lemma~4.1]{LS1}.
\end{proof}

\begin{proof}[Proof of Theorem~\ref{t:nonfill1}]
  Let $\xi _r$ be one of the contact structures on $M=M(g, 2g;
  (\alpha, 1))$ given by Figure~\ref{f:general}.  We shall argue as
  in~\cite[Theorem~1.1]{LS1}, therefore we shall need to find a
  $4$--manifold $Z=Z(g,2g; (\alpha , 1 ))$ with $b_2^+(Z)=0$,
  $\partial Z=-M$ and such that the intersection form $Q_Z$ does not
  embed into the diagonal lattice ${\mathbb {D}}_m=(\Z ^m , m(-1))$ for
  any $m$.
  
  We shall use a construction similar to the one given
  in~\cite[Proposition~4.4]{LS1}.  To this end, fix $g, d\in \N$ with
  $d(d+1)\leq 2g \leq (d+1)^2-2$, let $C\subset\CP^2$ be a smooth
  complex curve of degree $d+2$, and let $\widehat{\CP^2}$ be the
  blow--up of $\CP^2$ at $(d+2)^2-2g-1$ distinct points of $C$. Denote
  by $\widehat C\subset\widehat{\CP^2}$ the proper transform of $C$. Let
  ${\widetilde C}\subset\widehat{\CP^2}$ be a smooth, oriented surface
  obtained by adding $g - \frac 12 d(d+1)$ fake handles to $\widehat
  C$. Blow up $\widehat{\CP^2}$ at one more point of $\widetilde C$,
  then blow up repeatedly at distinct points of the last exceptional
  sphere until the corresponding proper transform in the resulting
  rational surface $X$ is an embedded sphere $S$ with
  self--intersection $-\alpha$. Define $Z$ as the complement in $X$ of
  a tubular neighborhood of $\widetilde C\cup S$. 
  
  The group $H_2(X;\Z)$ is generated by classes
  $h,e_1,e_2,\ldots,e_t$, where $h$ corresponds to the standard
  generator of $H_2(\CP^2;\Z)$ and the $e_i$'s are the classes of the
  exceptional curves. Let $q$ be a positive integer such that $2q\leq
  t$, and define $\La_q=(H_q,Q_q)$ as the intersection lattice given
  by the subgroup
\[
H_q=\langle e_1-e_2, e_2-e_3, \ldots , e_{2q-1}-e_{2q},
h-e_1-e_2-\ldots - e_q\rangle\subset H_2(X;\Z)
\] 
together with the restriction $Q_q$ of the intersection form
$Q_X$. 

As in the proof of~\cite[Proposition~4.4]{LS1}, the inequality $2g\leq
d(d+2)-1$ guarantees that $2(d+2)\leq t$, hence the 
lattice $\Lambda _{d+2}=(H_{d+2},
Q_{d+2})$ embeds into $(H_2(Z; \Z ), Q_Z)$. Since
by~\cite[Lemma~4.3]{LS1} $\Lambda _{d+2} =(H_{d+2}, Q_{d+2})$ does not
embed into any diagonal lattice ${\mathbb D}_m$, the same holds for
$(H_2(Z; \Z ), Q_Z)$.

By Corollary~\ref{c:weak}, a filling $(W,\om)$ would give rise to a negative
definite closed $4$--manifold $V=W\cup Z$ with nonstandard intersection form,
contradicting Donaldson's famous diagonalizability result \cite{Do1}.
\end{proof}

\section{Concluding remarks}
\label{s:nonfill}

With a little more work, essentially the same proof as the one given
in Section~\ref{s:homotopy} yields non--fillability for all structures
defined by Figure~\ref{f:general} on $M(g, n ; (\alpha , \beta ))$
and satisfying
\[
d(d+1)\leq 2g \leq n \leq d(d+2)-1
\]
for $g\geq 1$ and some integer $d$. In fact, a slightly more general 
argument in the computation of the $\Spin^c$ structures allows one 
to check that the statement of Theorem~\ref{t:moy} still holds.

\begin{figure}[ht]
\setlength{\unitlength}{1mm}
\begin{center}
\begin{picture}(100,55)
\put(12,8){$\frac 12(\al-r)$}
\put(71,8){$\frac 12(\al+r) - 1$}
\put(0,0){\includegraphics[width=10cm]{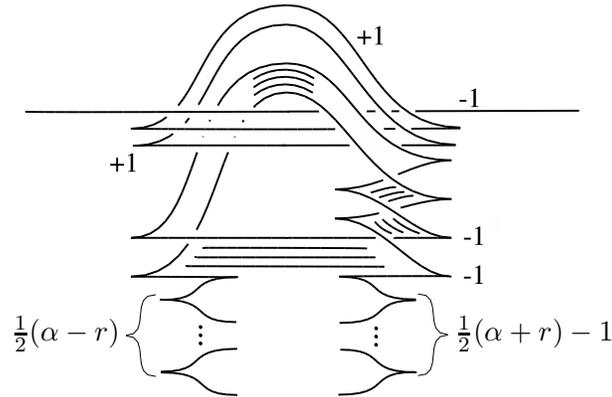}}
\end{picture}    
\end{center}
\caption{The contact structures $\xi^+_r$}
\label{f:final1} 
\end{figure}

\begin{figure}[ht]
\setlength{\unitlength}{1mm}
\begin{center}
\begin{picture}(100,55)
\put(8,9){$\frac 12(\al-r)-1$}
\put(72,9){$\frac 12(\al+r)$}
\put(0,0){\includegraphics[width=10cm]{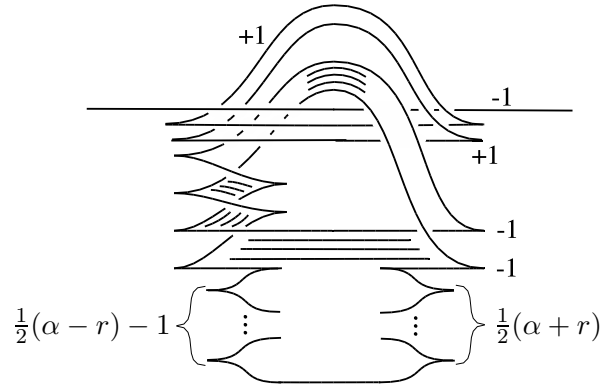}}
\end{picture}    
\end{center}
\caption{The contact structures $\xi^-_r$}
\label{f:final2} 
\end{figure}

In another direction, Theorem~\ref{t:nonfill1} generalizes to all
$M(g, n; (\alpha , 1))$ with $n \geq 2g>0$. In this case, one needs to
consider Figure~\ref{f:general} for $k=1$ and
\[
r_1=\frac{(n-2g+1)\al+1}{(n-2g+2)\al+1}.
\]
According to the algorithm described in Section~\ref{s:intro}, the
corresponding contact surgery can be expressed as a contact $(\pm
1)$--surgery by replacing the $r_1$--framed unknot $K$ with two
pushoffs of $K$, $n-2g$ pushoffs of a stabilization $K^\pm$ of $K$,
and one pushoff of $K^\pm$ stabilized $\al-1$ times. Depending on the
choice of stabilization of $K$, the result looks either like
Figure~\ref{f:final1} or Figure~\ref{f:final2}. Denoting by $r$
the rotation number of the last knot (after a choice of orientation),
this gives a contact structure $\xi^+_r$ for every $-\al<r\leq\al$ and
a contact structure $\xi^-_r$ for every $-\al\leq r < \al$ (and $r\equiv 
\alpha \bmod 2$ in both cases).

A computation as in Section~\ref{s:homotopy} gives
\[
\t_{\xi^\pm_r}=\t_{\rm can}+ \frac12 (r-\al-2\pm\al(n-2g)-
\al(n-2g))\PD(\mu).
\]
This already shows that the contact structures defined on $M(g,n;
(\alpha , 1))$ by Figure~\ref{f:general} are all distinct up to
homotopy, providing further evidence for Conjecture~\ref{c:differ}.

One can also compute the $3$--dimensional invariant $d_3([\xi^\pm_r])$
of the homotopy class $[\xi^\pm_r]$ of tangent $2$--plane fields
containing the contact structure $\xi^\pm_r$ (as discussed
in~\cite{LS2}), obtaining:
\[
d_3([\xi^{\pm}_r])=\frac{1}{4(n\alpha +1)}((n-2g)^2\alpha - r^2n\pm
2(n-2g)r) +\frac{2g-1}{2}.
\]

On the other hand, the statement of Theorem~\ref{t:moy} holds for all
contact structures defined on $M=M(g,n ; (\alpha , 1))$ by
Figure~\ref{f:general} for $n \geq 2g$. Therefore, the argument
of~\cite[Theorem~2.1]{L} and~\cite[Theorem~4.1]{LS2} applies, showing
that there is a unique homotopy class $\Xi(\t_{\xi^{\pm}_r})$ of
$2$--plane fields inducing the $\Spin^c$ structure $\t_{\xi^{\pm}_r}$
and which might potentially contain a fillable contact structure. The
proof of this observation rests on the fact that, assuming
Theorem~\ref{t:moy} to hold, the $3$--dimensional invariant of
$\Xi(\t_{\xi^{\pm}_r})$ is determined by some topological terms plus
an $\eta$--invariant of $(M, \t_{\xi^\pm_r})$ as follows.
  
By the formula preceding \cite[Section~3]{nic} 
(when $\rho (L)\neq
0$, which always holds in our case), the dimension $d_1$ of the
Seiberg--Witten moduli space with fixed boundary limit can be
expressed as 
\[
d_1=d_3(\Xi(\t_{\xi^\pm_r}))+\omega _{red}(\t_{\xi^\pm_r})-(2g-1), 
\]
where $d_3(\Xi(\t_{\xi^\pm_r}))$ is the $3$--dimensional invariant of
$\Xi(\t_{\xi^\pm_r})$ and $\omega _{red}(\t_{\xi^\pm_r})$ is given, in
the notations of~\cite{nic}, by the formula:
\begin{equation*}
\begin{split}
\frac{2g-1}{2}-\frac{l-sign(l)}{4}+ l\rho^\pm(1-\rho^\pm)-\rho^\pm
+\frac{1-\alpha}{2\alpha}(1-2\rho^\pm) 
+S(1,\alpha )\\+ F_{\rho^\pm}(\alpha, 1, \gamma)
+ 2S_{\rho^\pm} (1,\alpha,\gamma ).
\end{split}
\end{equation*}
In our situation we have:
\[
l=n+\frac{1}{\alpha },\quad 
{\rm sign}(l)=1, \quad 
\rho^\pm = \frac{\alpha(n\mp(n-2g))-r+1}{2n\alpha +2},
\]
\[
\gamma =\frac 12 (r+\alpha -2),\quad
S(1,\alpha )=\frac{\alpha ^2 +2}{12\alpha }-\frac 14,\quad
F_{\rho^\pm}(\alpha , 1, \gamma ) =\frac{\gamma +\rho^\pm}{\alpha},
\]
\[
S_{\rho^\pm}(1,\alpha,\gamma)= \frac{\alpha ^2-3\alpha (1+2\gamma )
+2(1+3\gamma +3 \gamma ^2)}{12\alpha }.
\]
This shows that 
\[
\omega _{red}(\t_{\xi^\pm_r})= -\frac{1}{4(n\alpha
+1)}((n-2g)^2\alpha - r^2n\pm 2(n-2g)r) + \frac{2g-1}{2}. 
\]

On the other hand, by the argument of~\cite[Theorem~2.1]{L} we have
\[
d_1=-1-b_1(M)=-1-2g,
\]
therefore
\[
d_3(\Xi(\t_{\xi^\pm_r}))=-\omega _{red}(\t_{\xi^\pm_r})-2,
\]
yielding
\[
d_3(\Xi(\t_{\xi^\pm_r}))
=\frac{1}{4(n\alpha +1)}((n-2g)^2\alpha - r^2n\pm 2(n-2g)r)
-\frac{2g+3}{2}.
\]
Since 
\[
d_3([\xi^{\pm}_r] )-d_3(\Xi(\t_{\xi^\pm_r}))=2g+1\neq 0, 
\]
none of the contact structures defined by Figure~\ref{f:general} on
$M(g, n; (\alpha , 1))$ ($n\geq 2g>0$) are symplectically
fillable.

We believe that the same idea should work for all the tight contact
structures given by Figure~\ref{f:general} (with the
constraints~\eqref{e:conditions}). The verification of
non--fillability, however, seems to be much more tedious in the
general case. The difficulty is number--theoretic in nature: it is
hard to see that $d_3([\xi])\neq d_3(\Xi(\t_\xi))$, because the
formulas involve sums which are hard to write in closed form.

\Addresses\recd

\begin{thebibliography}{AAA}

\bibitem{DG1}
{\bf F. Ding and H. Geiges},
{\it Symplectic fillability of tight contact structures on 
torus bundles}, \href{http://www.maths.warwick.ac.uk/agt/AGTVol1/agt-1-8.abs.html}{Algebr. Geom. Topol. {\bf1} (2001), 153--172.}

\bibitem{DG2} {\bf F. Ding and H. Geiges}, {\it A Legendrian surgery
    presentation of contact $3$--manifolds}, to appear in
    Proc. Cambridge Philos. Soc., \arxiv{math.SG/0107045}

\bibitem{DGS} {\bf F. Ding, H. Geiges and A. Stipsicz}, {\it Surgery
diagrams for contact $3$--manifolds}, to appear in Turkish J. Math.,
\arxiv{math.SG/0307237}

\bibitem{Do1}
{\bf S. K. Donaldson},
{\it An application of gauge theory to $4$--dimensional topology},
J. Differential Geom. {\bf 18} (1983), 279--315.

\bibitem{el}
{\bf Y. Eliashberg},
{\em Topological characterization of Stein manifolds of dimension $>2$},
International J. of Math. {\bf1} (1990), 29--46. 


\bibitem{EH}
{\bf J. Etnyre and K. Honda},
{\it Tight contact structures with no symplectic fillings},
Invent. Math. {\bf 148} (2002), no. 3, 609--626.


\bibitem{Gi} 
{\bf E.~Giroux}, 
{\it G\'eom\'etrie de contact de la dimension trois vers les dimensions
sup\'erior}, Proceedings of the ICM, Beijing 2002, vol. 2, 
405--414.


\bibitem{G}
{\bf R. Gompf},
{\it Handlebody constructions of Stein surfaces},
Ann. of Math. {\bf148} (1998), 619--693.

\bibitem{H2}
{\bf K. Honda},
{\it On the classification of tight contact structures, II.},
J. Differential Geom. {\bf55} (2000), 83--143.

\bibitem{PLpos}
{\bf P. Lisca},
{\it Symplectic fillings and positive scalar curvature},
Geom. Topol. {\bf2} (1998), 103--116.

\bibitem{L}
{\bf P. Lisca},
{\it On fillable contact structures up to homotopy},
Proc. Amer. Math. Soc. {\bf129} (2001), 3437--3444.

\bibitem{LS1}
{\bf P. Lisca and A. Stipsicz},
{\it An infinite family of tight, not semi--fillable 
contact three--manifolds}, 
\href{http://www.maths.warwick.ac.uk/gt/GTVol7/paper30.abs.html}%
{Geom. Topol. {\bf7} (2003), 1055--1073.}

\bibitem{LS2}
{\bf P. Lisca and A. Stipsicz},
{\it Tight, not semi--fillable contact circle bundles},
Math. Ann. {\bf328} (2004), 285--298.

\bibitem{MOY}
{\bf T. Mrowka, P. Ozsv\'ath, B. Yu},  
{\it Seiberg--Witten monopoles on Seifert fibered spaces}, 
Comm. Anal. Geom. {\bf 5} (1997), no. 4, 685--791. 

\bibitem{nic}
{\bf L. Nicolaescu},
{\it Finite energy Seiberg--Witten moduli spaces on $4$--manifolds bounding
Seifert fibrations}, Comm. Anal. Geom. {\bf8} (2000), 1027--1096.

\bibitem{orlik}
{\bf P. Orlik},
{\it Seifert Manifolds}, Lecture Notes in Mathematics, {\bf291}, 
Springer--Verlag (1972).

\bibitem{OSzF1} {\bf P. Ozsv\'ath and Z. Szab\'o}, {\it Holomorphic
    disks and topological invariants for closed three--manifolds}, to
  appear in Ann. Math., \arxiv{math.SG/0101206}
  
\bibitem{OSzF2} {\bf P. Ozsv\'ath and Z. Szab\'o}, {\it Holomorphic
    disks and three--manifold invariants: properties and
    applications}, to appear in Ann. Math., \arxiv{math.SG/0105202}
 
\bibitem{OSzF4} {\bf P. Ozsv\'ath and Z. Szab\'o}, {\it Holomorphic
    triangles and invariants of smooth $4$--manifolds},
  \arxiv{math.SG/0110169}

\bibitem{abs}
{\bf P. Ozsv\'ath and Z. Szab\'o},
{\it Absolutely graded Floer homologies and intersection forms
for four--manifolds with boundary},
Adv. Math. {\bf173} (2003), 179--261.

\bibitem{OSz6}
{\bf P. Ozsv\'ath and Z. Szab\'o},
{\it Heegaard Floer homologies and contact structures}, 
\arxiv{math.SG/0210127}

\bibitem{W}
{\bf A. Weinstein},
{\em Contact surgery and symplectic handlebodies},
Hokkaido Mathematical Journal {\bf20} (1991), 241--51.

\end{thebibliography}
\end{document}